\documentclass[a4paper]{article}
\usepackage{amsmath,amsthm,amssymb}
\usepackage[dvips]{graphicx}

\newtheorem{theorem}{Theorem} 

\newtheorem{proposition}[theorem]{Proposition}

\begin{document}
\title{Special Solutions of the Sixth Painlev\'e Equation \\
with Solvable Monodromy}
\author{Kazuo Kaneko \\
\small  Graduate School of Information Science and Technology, %\\
\small Osaka University\\[+12pt]
Shoji Okumura\\
\small  Graduate School of Science, Osaka University
}
\date{}
\maketitle
\begin{abstract}
We will study two types of special solutions of the sixth Painlev\'e 
equation, which are invariant under the symmetries obtained from the 
B\"acklund transformations. In most cases, the fixed points of the 
B\"acklund transformations are classical solution, but our solutions 
are not classical for generic parameters. We will calculate the linear 
monodromy of these solutions exactly, and we will characterize them on 
Fricke's cubic surface of monodromy. 
\end{abstract}
\section{Introduction}
R. Fuchs showed that the sixth Painlev\'e equation is represented as an 
isomonodromic deformation condition of the Fuchs type linear equation 
with four regular singular points~\cite{RF1}. Garnier showed that every 
type of the Painlev\'e equation is also obtained from monodromy preserving 
deformation of linear equation~\cite{Gar}. 
Nevertheless it is difficult to calculate the monodromy itself for the 
generic Painlev\'e functions. 

We call the monodromy data of the linear equation as a \textit{linear 
monodromy} of the Painlev\'e function. 
Our interests is in the Painlev\'e function whose corresponding linear 
monodromy can be determined exactly. In this paper, we call such 
Painlev\'e functions \textit{monodromy solvable}.

One example of monodromy solvable Painlev\'e functions is Umemura's 
classical solutions~\cite{HU1}, which are fixed points of the B\"acklund 
transformations. For any classical solutions, we can calculate the 
linear monodromy. 
But conversely, there exist some non-classical Painlev\'e functions 
whose linear monodromy can be calculated. 
We are interested in non-classical monodromy solvable solutions. 

R. Fuchs is the first to find a non-classical monodromy solvable 
solution~\cite{RF2}. 
He calculates the linear monodromy of Picard's solution \cite{EP}, 
which satisfies the sixth Painlev\'e equation with a special parameter. 
His work is found again in \cite{MMA} recently. 

Another example of monodromy solvable solution is a symmetric solution. 
The first, second and fourth Painlev\'e equations 
\begin{alignat*}{2}
&(P1)&\qquad &\frac{d^2 y}{dt^2}=6y^2+t,\\
&(P2)&\qquad &\frac{d^2 y}{dt^2}=2y^{3}+ty+\alpha,\\
&(P4)&\qquad &\frac{d^2 
y}{dt^2}=\frac{1}{2y}\left(\frac{dy}{dt}\right)^{2}+\frac{3}{2} 
 y^3+4ty^{2}+2(t^{2}-\alpha)y+\frac\beta y
\end{alignat*}
are invariant under the transformations
\begin{alignat*}{4}
&(P1)&\qquad  &t\to \rho t,&\quad y&\to \rho^3 y& \quad(\rho^5=1), \\
&(P2)&\qquad  &t\to \omega t,&\quad y&\to \omega^2 y&\quad 
(\omega^3=1), \\ 
&(P4)&\qquad  &t \to -t,&\quad y&\to -y.
\end{alignat*}
There exist symmetric solutions which are invariant under the action of 
the cyclic groups. The symmetric solutions are studied by 
Kitaev~\cite{AVK} for (P1) and (P2), and by Kaneko~\cite{KK1} for (P4). 
For these cases, $(y,t)=(0,0)$ is a fixed point of the transformations, 
and there exist symmetric solutions with initial values at 
the fixed point $y(0)=0$. Kitaev calculated the linear monodromy of the 
symmetric solutions of (P1) and (P2) by adopting the fixed point as the 
initial condition. In the same way, Kaneko calculated the linear 
monodromy of the symmetric solution of (P4). We remark that these 
symmetric solutions are not classical for generic parameters. 

In this paper, we will construct monodromy solvable solutions of the 
sixth Painlev\'e equation. 

For generic parameter, we cannot find a simple symmetry as above. 
However, for some restricted parameters, there exist two types of 
symmetries obtained from the B\"acklund transformations. 
One is 
\[
\sigma_{1}:\ t\to 1-t,\qquad y\to 1-y,\qquad \beta \to -\gamma,\qquad \gamma \to -\beta,
\]
and the other is 
\[
 \sigma_{2}\circ\sigma_{1}:\ t\to 1/(1-t),\qquad y\to 1/(1-y),\qquad
 \alpha\to-\beta,\qquad
 \beta\to-\gamma,\qquad
 \gamma\to\alpha. 
\]
The sixth Painlev\'e equation is invariant under the action $\sigma_{1}$ 
if $-\beta=\gamma$, and is invariant under the action 
$\sigma_{2}\circ\sigma_{1}$ if $\alpha=-\beta=\gamma$. 

In the section 3, we will show there exist symmetric solutions, which 
are fixed points of these actions. 
In many cases, the fixed points of the B\"acklund transformations 
are classical solutions. However, our symmetric solutions are not 
classical for generic parameters.

In the section 4, we will show that the linearizations of these symmetric 
solutions are reduced to the Gauss hypergeometric equations. In the 
section 5, we will calculate the linear monodromy of our symmetric solutions 
explicitly. 

If the monodromy matrices $\left\{M_{0},M_{t},M_{1},M_{\infty}\right\}$
are elements of $SL(2,\mathbb{C})$, then
$$
p_{j}=\mathrm{tr}M_{j}, \qquad p_{jk}=\mathrm{tr} M_{j}M_{k},
$$
satisfies the following relation:
\begin{multline}
p_{01}p_{1t}p_{t0}+p_{01}^2+p_{1t}^2+p_{t0}^2
-(p_{0}p_{1}+p_{t}p_{\infty})
p_{01}-(p_{1}p_{t}+p_{0}p_{\infty})p_{1t}-(p_{t}p_{0}+p_{1}p_{\infty})p_{t0}
\\
+p_{0}^2+p_{1}^2+p_{t}^2+p_{\infty}^{2}
+p_{0}p_{t}p_{1}p_{\infty}-4=0.\label{eq-mon}
\end{multline}
It is known that the monodromy matrices 
$\{M_{0},M_{t},M_{1},M_{\infty}\}$ are  
determined by 
$\{p_{0}, p_{1}, p_{t}, p_{\infty}; p_{01}, p_{1t}, p_{t0}\}$ 
upto gauge transformations(\cite{Boa}, \cite{Iwa}, \cite{Jim}, \cite{Fri}). 
%Therefore the monodromy is on the cubic surface \eqref{eq-mon}. 
We call \eqref{eq-mon} as Fricke's cubic surface of monodromy. 
In the section 6, we will 
characterize the monodromy of our symmetric solutions on Fricke's cubic 
surface of monodromy. 

Andreev and Kitaev constructed many algebraic solutions of 
the sixth Painlev\'e equation, by rational transformations of the 
hypergeometric equations~\cite{AK}. By the rational transformations of 
confluent hypergeometric equation, Ohyama and Okumura constructed 
algebraic solutions and symmetric solutions of the 
first to the fifth Painlev\'e equations~\cite{Oh}. 
We remark that Boach also constructed many algebraic solutions of the 
sixth Painlev\'e equation, whose linear monodromy is isomorphic to the 
complex reflection group~\cite{Boa}. 

The solutions of this paper essentially appeared 
in~\cite{AK}. But, since their interest are in algebraic solutions, they 
did not go to detail of these solutions. In this paper, we will 
calculate the linear monodromy explicitly, and we will show that they 
are symmetric solutions. 

\bigskip

The authors give thanks to Professor Yousuke Ohyama 
for fruitful discussions and suggestions. 

\section{Linear Problem}\label{LP}
The sixth Painlev\'e equation \eqref{p6} is obtained from the 
isomonodromic deformation equation of 
linear equation with four regular singular points. We will use 
Garnier-Okamoto's linearization \cite{KO}. 
\begin{align}
&\frac{\partial^2 \psi}{\partial x^2}+p(x,t)\frac{\partial \psi}{\partial 
x}+q(x,t)\psi=0, \label{OG-1}\\
&\frac{\partial\psi}{\partial t}=a(x,t)\frac{\partial \psi}{\partial x}+
b(x,t)\psi, \label{OG-2}
\end{align}
where
\begin{align}
p(x,t)&=\frac{1-\alpha_{4}}{x}+\frac{1-\alpha_{3}}{x-1}+\frac{1-\alpha_{0}}{x-t}-\frac{1}{x-y},\\
q(x,t)&=\frac{\alpha_{2}(\alpha_{1}+\alpha_{2})}{x(x-1)}-
\frac{t(t-1)H_{VI}}{x(x-1)(x-t)}+\frac{y(y-1)z}{x(x-1)(x-y)},\\
a(x,t)&=\frac{y-t}{t(t-1)}\frac{x(x-1)}{x-y},\\
b(x,t)&=\frac{(1-\alpha_{4}-\alpha_{3}-\alpha_{0})(y-t)}{2t(t-1)}-
\frac{y(y-1)(y-t)z}{t(t-1)(x-y)},
\end{align}
and
\begin{multline}
 H_{VI}=\frac 1{t(t-1)}\Bigl[
y(y-1)(y-t)z^2-\Bigl\{\alpha_4(y-1)(y-t)\\{}+\alpha_3y(y-t)+(\alpha_0-1)y(y-1)
\Bigr\}z 
+\alpha_2(\alpha_{1}+\alpha_{2})(y-t)
\Bigr].
\end{multline}
The Riemann scheme of \eqref{OG-1} is
\[\left\{x;
 \begin{matrix}
	0	& 1	& t	& y	& \infty \\
	0	& 0	& 0	& 0	& \alpha_{2}\\
	\alpha_{4}&\alpha_{3}&\alpha_{0}&2&\alpha_{1}+\alpha_{2}
 \end{matrix}\right\}.
\]
From the compatibility condition of \eqref{OG-1} and \eqref{OG-2},
we have the following Hamiltonian system
\[
 \frac{dy}{dt}=\frac{\partial H_{VI}}{\partial z}, \qquad
 \frac{dz}{dt}=-\frac{\partial H_{VI}}{\partial y},
\]
that is
\begin{align}
\begin{aligned}
t(t-1)\frac{dy}{dt}=&2y(y-1)(y-t)z- 
 (\alpha_{0}-1)y(y-1)- 
 \alpha_3y(y-t)- 
 \alpha_4(y-1)(y-t),\\[+3pt]
t(t-1)\frac{dz}{dt}=&-\Bigl(y(y-1)+(y-1)(y-t)+y(y-t)\Bigr)z^2 \\
&\qquad\qquad
+\Bigl((2y-1)(\alpha_{0}-1)+(2y-t)\alpha_{3}+(2y-t-1)\alpha_{4}\Bigr)z-
  (\alpha_1 + \alpha_2)\alpha_2 .
\end{aligned}\label{P6H}
\end{align}
Let us eliminate $z$ from this system, and we have the sixth Painlev\'e 
equation 
\begin{multline}
   \frac{d^2 y}{d t^2} = \frac{1}{2}\left( \frac{1}{y} +
       \frac{1}{y-1} + \frac{1}{y-t} \right) \left(\frac{d y}{d t}\right)^2
            - \left( \frac{1}{t}+\frac{1}{t-1}
               +\frac{1}{y-t} \right)\frac{d y}{d t}  \\
      + \frac{y\left(y-1\right)\left(y-t\right)}{t^2\left(t-1\right)^2}
        \left\{ \alpha + \beta\frac{t}{y^2} +
            \gamma\frac{t-1}{\left(y-1\right)^2 }+
             \delta \frac{t\left(t-1\right)}{\left(y-t\right)^2} \right\},
\label{p6}
\end{multline}
where
\[
 \alpha=\frac{\alpha_{1}^{2}}{2},\qquad
 \beta=-\frac{\alpha_{4}^{2}}{2},\qquad
 \gamma=\frac{\alpha_{3}^{2}}{2},\qquad
 \delta=\frac{1-\alpha_{0}^{2}}{2}.
\]

\bigskip

We will calculate the linear monodromy of some special Painlev\'e 
functions. 

\section{B\"acklund Transformations of the Sixth Painlev\'e Equation}
\label{sym}
For the sixth Painlev\'e equation,
we cannot find a simple symmetry for general parameters. However, for 
some restricted parameters, there exist some symmetries obtained from 
the B\"acklund transformations. And there exist symmetric solutions 
of the B\"acklund transformations. 

\bigskip

The sixth Painlev\'e equation \eqref{P6H} is invariant under the 
following transformations\cite{NY}:
\begin{center}
\begin{tabular}{|c||ccccc|cc|c|}\hline
{} &$ \alpha_0 $&$ \alpha_1 $&$ \alpha_2 $&$ \alpha_3 $&$ \alpha_4 $&$ y $&$ z $&$ t$\\
\hline
$s_0 $&$ -\alpha_0 $&$ \alpha_1 $&$ \alpha_0+\alpha_2 $&$ \alpha_3 $&$ \alpha_4 $&$ y $&$ z-\frac{\alpha_0}{y-t} $&$ t$\\$
s_1 $&$ \alpha_0 $&$ -\alpha_1 $&$ \alpha_1+\alpha_2 $&$ \alpha_3 $&$ \alpha_4 $&$ y $&$ z$&$ t$\\$s_2 $&$ \alpha_0+\alpha_2 $&$ \alpha_1+\alpha_2 $&$ -\alpha_2
$&$ \alpha_2+\alpha_3 $&$ \alpha_2+\alpha_4 $&$ y+\frac{\alpha_2}{z} $&$ z $&$ t$\\$
s_3 $&$ \alpha_0 $&$ \alpha_1 $&$ \alpha_2+\alpha_3 $&$ -\alpha_3 $&$ \alpha_4 $&$ y
$&$ z-\frac{\alpha_3}{y-1} $&$ t$\\$s_4 $&$ \alpha_0 $&$ \alpha_1 $&$ \alpha_2+\alpha_4
$&$ \alpha_3 $&$ -\alpha_4 $&$ y $&$ z-\frac{\alpha_4}{y} $&$ t$\\
\hline
$ \pi_1 $&$ \alpha_3$&$ \alpha_4 $&$ \alpha_2 $&$ \alpha_0 $&$ \alpha_1 $&$ \frac{t}{y}
$&$ -\frac{y\left(yz+\alpha_2\right)}{t} $&$ t$\\
$\pi_2 $&$ \alpha_1 $&$ \alpha_0 $&$ \alpha_2 $&$ \alpha_4 $&$ \alpha_3
$&$ \frac{t\left(y-1\right)}{y-t}
$&$ 
-\frac{\left(y-t\right)^2z+\left(y-t\right)\alpha_2}{\left(t-1\right)
t}$&$ t$\\ 
\hline
$\sigma_1 $&$ \alpha_0 $&$ \alpha_1 $&$ \alpha_2 $&$ \alpha_4 $&$ \alpha_3 $&$ 1-y $&$ -z $&$ 1-t$\\$
\sigma_2 $&$ \alpha_0 $&$ \alpha_4 $&$ \alpha_2 $&$ \alpha_3 $&$ \alpha_1 $&$ \frac{1}{y} $&$ -y\left(yz+\alpha_2\right) $&$ \frac{1}{t}$\\$
\sigma_3 $&$ \alpha_4 $&$ \alpha_1 $&$ \alpha_2 $&$ \alpha_3 $&$ \alpha_0 $&$ \frac{t-y}{t-1} $&$ \left(1-t\right)z $&$ \frac{t}{t-1}$\\
\hline
\end{tabular}\end{center}
where $\alpha_{0},\alpha_{1},\alpha_{2},\alpha_{3}$ and $\alpha_{4}$ are the 
parameters in \eqref{P6H}, and satisfy the following relation:  
\[
 \alpha_{0}+\alpha_{1}+2\alpha_{2}+\alpha_{3}+\alpha_{4}=1.
\]
We remark that 
$$\langle s_{0}, s_{1}, s_{2}, s_{3}, s_{4}, \pi_{1}, \pi_{2}\rangle$$
is an extended affine Weyl group of type $D_{4}^{(1)}$, and 
$$\langle s_{0}, s_{1}, s_{2}, s_{3}, s_{4}, \pi_{1}, \pi_{2}, 
\sigma_{1},\sigma_{2},\sigma_{3}\rangle$$ is an extended affine Weyl 
group of type $F_{4}^{(1)}$.

\bigskip

If $\alpha_{3}=\alpha_{4}$, \eqref{P6H} is invariant under 
the action of $\sigma_{1}$, and if $\alpha_{1}=\alpha_{3}=\alpha_{4}$, 
\eqref{P6H} is invariant under the action of $\sigma_{2}\circ\sigma_{1}$. 
In the following, we will study the symmetric solutions of $\sigma_{1}$ 
and $\sigma_{2}\circ\sigma_{1}$.

\subsection{Symmetric Solution of $\sigma_{1}$}

For the B\"acklund transformation $\sigma_{1}$: 
\[
 y \to 1-y, \qquad z \to -z, \qquad t \to 1-t,
\]
the sixth Painlev\'e equation is invariant with the condition 
$\alpha_{3}=\alpha_{4}$, and one of the fixed point is
\[
 t=1/2, \quad y=1/2, \quad z=0.
\]
Since this is not a critical point of $\eqref{P6H}$, 
there exist a unique holomorphic solution of \eqref{P6H} with the 
initial condition 
\[
 y(1/2)=1/2, \qquad z(1/2)=0.
\]
And if $\alpha_{3}=\alpha_{4}$, this solution is expanded as follows: 
\begin{align}
y(t)=\sum_{k=0}^\infty a_{k}\left(t-\tfrac 12\right)^{k},\qquad
z(t)=\sum_{k=0}^{\infty} b_{k}\left(t-\tfrac 12\right)^{k},
\label{sym2}
\end{align}
where
\begin{align*}
&a_{0}=\frac 12,\quad a_{1}=1-\alpha_{0},\quad a_{2}=0,\\
&\qquad a_{3}=-\frac 43 \alpha_{0}\bigl\{
2\alpha_{2}(\alpha_{1}+\alpha_{2})-(1-\alpha_{0})(2\alpha_{2}+\alpha_{1}+1)
\bigr\},\quad a_{4}=0, \quad \cdots,\\[+2ex]
&b_{0}=0,\quad b_{1}=4\alpha_{2}(\alpha_{1}+\alpha_{2}),\quad b_{2}=0,\\
&\qquad
b_{3}=-\frac{16}3 
\alpha_{2}(\alpha_{1}+\alpha_{2})\bigl\{
\alpha_{2}(\alpha_{1}+\alpha_{2})-(\alpha_{1}+2\alpha_{2})(1-2\alpha_{0})
+\alpha_{0}-2 \bigr\},\quad b_{4}=0, \quad
\cdots.
\end{align*}

\begin{proposition}
The Painlev\'e function \eqref{sym2} is a symmetric solution of 
$\sigma_{1}$:
\[
 y(1-t)=1-y(t),\qquad
 z(1-t)=-z(t).
\]
\end{proposition}
\begin{proof}
Let us set
\[
 \tau=t-1/2, \qquad \lambda=y-1/2, \qquad \mu=z,
\]
then the action of $\sigma_{1}$ maps $\tau\to -\tau$, $\lambda\to 
-\lambda$, and $\mu\to -\mu$. 
\eqref{P6H} becomes the following Hamiltonian system 
$\frac{d\lambda}{d\tau}=\frac{\partial  
K}{\partial \mu}$, $\frac{d\mu}{d\tau}=-\frac{\partial K}{\partial 
\lambda}$: 
\begin{align}
\begin{aligned}
\left(\tau^2-\tfrac{1}{4}\right)\frac{d \lambda}{d \tau}&=
	2\left(\lambda^2-\tfrac{1}{4}\right)\mu\left(\lambda-\tau\right)
	-\left(\lambda^2-\tfrac{1}{4}\right)\left(\alpha_{0}-1\right)
	-2\lambda\left(\lambda-\tau\right)\alpha_{3},\\
\left(\tau^2-\tfrac{1}{4}\right)\frac{d \mu}{d \tau}&=
	-\left(\left(\lambda^2-\tfrac{1}{4}\right)
	+2\lambda\left(\lambda-\tau\right)\right)\mu^2
	+2\left(\lambda(\alpha_{0}-1)+(2\lambda-\tau)\alpha_{3}
	\right)\mu
	-\alpha_{2}(\alpha_{1}+\alpha_{2}).
\end{aligned}\label{H2mod}
\end{align}
The Hamiltonian $K$ is
\[
 \left(\tau^2-\tfrac 1 4\right)K=
\mu^2(\lambda-\tau)\left(\lambda^2-\tfrac 14\right)-\mu\left\{
\left(\lambda^2-\tfrac 1 4\right) (\alpha_{0}-1)+2\lambda(\lambda-\tau)\alpha_{3}\right\}+(\lambda-
\tau)\alpha_{2}(\alpha_{1}+\alpha_{2}) .
\]
Then \eqref{sym2} corresponds to the following solution of \eqref{H2mod}: 
\begin{align*}
\lambda = (1-\alpha_0)\tau + O(\tau)^3,\qquad
\mu     =4\alpha_2(\alpha_1+\alpha_2)\tau + O(\tau)^3.
\end{align*}
Higher order expansions are determined from 
\eqref{H2mod} inductively, and it has only odd powers of $\tau$. 
Therefore it determines symmetric solution. Thus \eqref{sym2} is a 
symmetric solution. 
\end{proof}

By the B\"acklund transformations $s_{0}$, $s_{1}$, $s_{2}$ and 
$\pi_{2}$, which are commutative with $\sigma_{1}$, we have four 
symmetric solutions around $t=1/2$. 

\begin{proposition} \label{sym2-prop}
If $\alpha_3=\alpha_4$, 
all of the symmetric solutions on $\sigma_{1}$ around $t=1/2$ are
\begin{align}
&
\left\{
\begin{aligned}
y- \tfrac 12 &= (1-\alpha_0)\left(t-\tfrac 12\right) 
	+ O\left(t-\tfrac 12\right)^3,
\\
z&=4\alpha_2(\alpha_1+\alpha_2)\left(t-\tfrac 12\right)
	+ O\left(t-\tfrac 12\right)^3,
\end{aligned}\right.\tag{S2-1}\label{sym2-1}
\\
&
\left\{
\begin{aligned} 
y-\tfrac 12&= (1+\alpha_0)\left(t-\tfrac 12\right) 
	+ O\left(t-\tfrac 12\right)^3,
\\
z&=\left(t-\tfrac 12\right)^{-1}
	+ O\left(t-\tfrac 12\right),
\end{aligned}\right.\tag{S2-2}\label{sym2-2}
\\ 
&\left\{
\begin{aligned}
y-\tfrac 12&= \tfrac 1{4\alpha_1}\left(t-\tfrac 12\right)^{-1} 
	+O\left(t-\tfrac 12\right),
\\
z&=-4\alpha_1\alpha_2\left(t-\tfrac 12\right)
	+ O\left(t-\tfrac 12\right)^3,
\end{aligned}\right.\tag{S2-3}
\label{sym2-3}
\\
&\left\{
\begin{aligned}
y-\tfrac 12&= \tfrac {-1}{4\alpha_1}\left(t-\tfrac 12\right)^{-1} 
	+ O\left(t-\tfrac 12\right) ,
\\
z&=4\alpha_1(\alpha_1+\alpha_2)\left(t-\tfrac 12\right)
	+ O\left(t-\tfrac 12\right)^3 ,
\end{aligned}\right.\tag{S2-4}
\label{sym2-4}
\end{align}
where \eqref{sym2-1} is \eqref{sym2}.

The actions of $s_{0}$, $s_{1}$, $s_{2}$ and $\pi_{2}$ interchange the 
symmetric solutions as follows:
\begin{center}
\begin{picture}(130,85)
 \put(0,20){\makebox(10,10){$\scriptstyle s_1$}}
 \put(15,25){\oval(8,8)[l]}
 \put(15,21){\vector(1,0){4}}
 \put(15,29){\vector(1,0){4}}
 \put(20,20){\makebox(30,10){\eqref{sym2-2}}}
 \put(50,20){\makebox(30,10){$\scriptstyle s_2$}}
 \put(55,25){\oval(8,8)[r]}
 \put(55,21){\vector(-1,0){4}}
 \put(55,29){\vector(-1,0){4}}
 \put(80,20){\makebox(30,10){\eqref{sym2-4}}}
 \put(75,25){\oval(8,8)[l]}
 \put(75,21){\vector(1,0){4}}
 \put(75,29){\vector(1,0){4}}
 \put(115,25){\oval(8,8)[r]}
 \put(115,21){\vector(-1,0){4}}
 \put(115,29){\vector(-1,0){4}}
 \put(110,20){\makebox(30,10){$\scriptstyle s_0$}}
 \put(35,10){\vector(0,1){9}}
 \put(95,10){\vector(0,1){9}}
 \put(35,10){\line(1,0){60}}
 \put(60,0){\makebox(10,10){$\scriptstyle \pi_{2}$}}
 \put(35,36){\vector(0,1){28}}
 \put(35,64){\vector(0,-1){28}}
 \put(36,45){\makebox(10,10){$\scriptstyle s_{0}$}}
 \put(95,35){\vector(0,1){30}}
 \put(95,65){\vector(0,-1){30}}
 \put(96,45){\makebox(10,10){$\scriptstyle s_{1}$}}
 \put(20,65){\makebox(30,10){\eqref{sym2-1}}}
 \put(80,65){\makebox(30,10){\eqref{sym2-3}}}
 \put(50,72){\vector(1,0){30}}
 \put(50,68){\vector(1,0){30}}
 \put(80,72){\vector(-1,0){30}}
 \put(80,68){\vector(-1,0){30}}
 \put(60,74){\makebox(10,10)[b]{$\scriptstyle s_{2}$}}
 \put(60,55){\makebox(10,10)[t]{$\scriptstyle \pi_{2}$}}
 \put(0,65){\makebox(10,10){$\scriptstyle s_1$}}
 \put(15,70){\oval(8,8)[l]}
 \put(15,66){\vector(1,0){4}}
 \put(15,74){\vector(1,0){4}}
 \put(115,70){\oval(8,8)[r]}
 \put(115,66){\vector(-1,0){4}}
 \put(115,74){\vector(-1,0){4}}
 \put(110,65){\makebox(30,10){$\scriptstyle s_0$}}
\end{picture}
\end{center}
\end{proposition}

In section \ref{MS2}, we will show these symmetric solutions
are monodromy solvable. 

\subsection{Symmetric Solution of $\sigma_{2}\circ \sigma_{1}$}
For the B\"acklund transformation $\sigma_{2}\circ\sigma_{1}$: 
\[
 y\to\frac 1{1-y},\qquad
 z\to-(1-y)(-z(1-y)+\alpha_{2}),\qquad
 t\to\frac 1{1-t},
\]
the sixth Painlev\'e equation is invariant with the condition 
$\alpha_{1}=\alpha_{3}=\alpha_{4}$, and one of the fixed point is
\[
 y=-\omega^2,\qquad
 z=\frac{2\omega+1}3\alpha_{2},\qquad
 t=-\omega^2.
\]
Since this is not a critical point of $\eqref{P6H}$, 
there exist a unique holomorphic solution of \eqref{P6H} with the initial 
condition 
\[
 y(-\omega^2)=-\omega^{2}, \qquad z(-\omega^2)=\frac{2\omega+1}3\alpha_{2}.
\]
And if $\alpha_{1}=\alpha_{3}=\alpha_{4}$, this solution is expanded as 
follows: 
\begin{align}
\begin{aligned}
&y(t)=-\omega^{2}+(1-\alpha_{0})(t+\omega^2)+\frac{1+2\omega}3\alpha_{0}(1-\alpha_{0})(t+\omega^2)^2+O(t+\omega^2)^3,\\
&z(t)=\frac{2\omega+1}3 \alpha_{2}-\frac{\alpha_{2}}3(1-\alpha_{0})
(t+\omega^2)+O(t+\omega^2)^2.
\end{aligned}\label{sym3}
\end{align}

\begin{proposition}
The Painlev\'e function \eqref{sym3} is a symmetric solution of 
$\sigma_{2}\circ\sigma_{1}$: 
\[
 y\left(\frac 1{1-t}\right)=\frac 1{1-y(t)},\qquad
 z\left(\frac 1{1-t}\right)=-(1-y(t))\left\{-z(t)(1-y(t))+\alpha_{2}\right\}.
\]
\end{proposition}
\begin{proof}
To prove the symmetry of $\sigma_{2}\circ\sigma_{1}$, we set
\begin{align*}
\tau=\frac{-\omega t-1}{t+\omega},\qquad
\lambda=\frac{-\omega y-1}{y+\omega},
\end{align*}
which diagonalize the action of 
$\sigma_2 \circ \sigma_1$: $\tau\to\omega\tau$, 
$\lambda\to\omega\lambda$. 
Then the sixth Painlev\'e equation \eqref{p6} becomes
\begin{multline}
(\lambda-\tau)\frac{d^2\lambda}{d\tau^2}=\frac{4\lambda^3-3\tau\lambda^2+1}
{2(\lambda^3+1)}\left(\frac{d\lambda}{d\tau}\right)^2+
\frac{-3\tau^2\lambda+2\tau^3-1}{\tau^3+1}
\frac{d\lambda}{d\tau}\\
{}+\frac{(\lambda^3+1)^2(\tau^3+1)(1-\alpha_0^2)+9(\lambda-\tau)^2
(\lambda^4-2\tau\lambda^3-2\lambda+\tau)\alpha_1^2}
{2(\lambda^3+1)(\tau^3+1)^2}, \label{P6 modified}
\end{multline}
which is equivalent to the following Hamiltonian system
$\frac{d\lambda}{d\tau}=\frac{\partial K}{\partial \mu}$,
$\frac{d\mu}{d\tau}=-\frac{\partial K}{\partial \lambda}$:
\begin{align}
\begin{aligned}
 (\tau^{3}+1)\frac{d\lambda}{d\tau}=&
2\left( 1 + \lambda^3 \right) \mu \left( \lambda  - \tau  \right)
+3\lambda^2\left( \lambda  - \tau  \right)\alpha_1
+\left( 1 + \lambda^3 \right) \left( 1 + \alpha_0 \right)
,\\
 (\tau^{3}+1)\frac{d\mu}{d\tau}=&
\mu^2\left( - 4\lambda^3 + 3\lambda^2\tau -1 \right)
-3\lambda \mu \left( \lambda \left( 1 + \alpha_0 \right)  + 
    \left( 3\lambda  - 2\tau  \right) \alpha_1 \right)
- \lambda {\left(  \alpha_2 -1\right) }^2 \\[-6pt]
&\hspace{5cm}+\frac{ \alpha_2-1 }{2}\left( \left( \lambda  + \tau  \right) \left( 1 + \alpha_0 \right)  + 
      3\left( \lambda  - \tau  \right) \alpha_1 \right)  .
\end{aligned}\label{H6 modified}
\end{align}
The Hamiltonian $K$ is
\begin{multline*}
 (\tau^{3}+1)K=\mu^2(\lambda^{3}+1)(\lambda-\tau)+
	\mu\Bigl\{
		(\lambda^3+1)(1+\alpha_{0})+3\lambda^2(\lambda-\tau)\alpha_{1}
\Bigr\}\\ 
	+\frac{1+\alpha_{0}+3\alpha_{1}}4
\lambda\Bigl\{
	(\lambda+\tau)(1+\alpha_{0})+3(\lambda-\tau)\alpha_{1}
\Bigr\},
\end{multline*}
and $\mu$ is determined by 
\begin{multline*}
\mu=\frac{t(t-1)(y+\omega)^2}{2y(y-1)(y-t)(1-\omega^2)}\frac{dy}{dt}
+\frac{(\omega-1)(t+\omega)^2(\alpha_{0}+1)}{6(y-t)}+
\frac{(1-\omega^2)\alpha_{1}}{6y}+
\frac{\omega^2(1-\omega^2)\alpha_{1}}{6(y-1)}\\
+
\frac{(1-\omega)\alpha_{1}y}{6}
+\frac 1 6 \left\{
(\omega-1)(t+\omega)(1+\alpha_{0})+(\omega^2-1)\alpha_{1}
\right\}.
\end{multline*}

By the action of $\sigma_{2}\circ\sigma_{1}$, we have 
\[
 \tau\to\omega\tau, \qquad
 \lambda\to\omega\lambda,\qquad
 \mu\to\omega^2\mu,\qquad
 K\to\omega^2 K.
\]
There exists a solution of \eqref{H6 modified}:
\begin{align}
\begin{aligned}
\lambda&=\left( 1 - \alpha_0 \right) \tau  %+ 
+O(\tau)^4,
\\
\mu&=\tau^{-1}+O(\tau)^2, 
\end{aligned}\label{sym3mod}
\end{align}
which are corresponding to \eqref{sym3}. And inductively, we have 
\eqref{sym3mod} in the following forms 
\begin{align*}
&\lambda=\sum a_{n} \tau^{3n+1}, 
&\mu=\sum b_{n} \tau^{3n-1}.
\end{align*}
Therefore it is a symmetric solution. And thus \eqref{sym3} is a symmetric 
solution. 
\end{proof}

By the B\"acklund transformations $s_{0}$ and $s_{2}$, which are 
commutative with $\sigma_{2}\circ\sigma_{1}$, we have three symmetric 
solutions around $t=-\omega^2$. 

\begin{proposition}
If $\alpha_{1}=\alpha_{3}=\alpha_{4}$, 
all of symmetric solutions on $\sigma_{2}\circ\sigma_{1}$ around 
$t=-\omega^2$ are 
\begin{align}
&\left\{
\begin{aligned}
y&=-\omega+O(t+\omega^2)^{2},\\
z&=-\tfrac 1 3 (1+2\omega)\alpha_{2}
	-\tfrac 1 6 (\alpha_{0}-3\alpha_{1}-1)\alpha_{2}(t+\omega^2)
	+O(t+\omega^2)^{2},
\end{aligned}\right.\tag{S3-1}\label{sym3-1}\\
&\left\{
\begin{aligned} 
y&=-\omega^2+(1-\alpha_{0})(t+\omega^2)+O(t+\omega^2)^2,\\
z&=\tfrac 1 3 (1+2\omega)\alpha_{2}+
	\tfrac 1 3(\alpha_{0}-1)\alpha_{2}(t+\omega^2)+O(t+\omega^2)^{2},
\end{aligned}\right.\tag{S3-2}\label{sym3-2}\\
&\left\{
\begin{aligned}
y&=-\omega^2+(1+\alpha_{0})(t+\omega^2)+O(t+\omega^2)^{2},\\
z&=(t+\omega^2)^{-1}+\tfrac 1 2 (1+2\omega)(\alpha_{0}-\alpha_{1}+1)
	+O(t+\omega^2),
\end{aligned}\right.\tag{S3-3}\label{sym3-3}
\end{align}
where \eqref{sym3-2} is \eqref{sym3}. 
The actions of $s_{0}$ and $s_{2}$ interchange the symmetric solutions 
as follows: 
\begin{center}
\begin{picture}(190,15)
 \put(0,0){\makebox(10,10){$\scriptstyle s_0$}}
 \put(15,5){\oval(8,8)[l]}
 \put(15,1){\vector(1,0){4}}
 \put(15,9){\vector(1,0){4}}
 \put(20,0){\makebox(30,10){\eqref{sym3-1}}}
 \put(50,5){\makebox(30,10){$\scriptstyle s_2$}}
 \put(50,5){\vector(1,0){30}}
 \put(80,5){\vector(-1,0){30}}
 \put(80,0){\makebox(30,10){\eqref{sym3-2}}}
 \put(110,5){\makebox(30,10){$\scriptstyle s_0$}}
 \put(110,5){\vector(1,0){30}}
 \put(140,5){\vector(-1,0){30}}
 \put(140,0){\makebox(30,10){\eqref{sym3-3}}}
 \put(175,5){\oval(8,8)[r]}
 \put(175,0){\vector(-1,0){4}}
 \put(175,9){\vector(-1,0){4}}
 \put(180,0){\makebox(10,10){$\scriptstyle s_2$}}
\end{picture}
\end{center}
There are more three symmetric solutions around $t=-\omega$:
\begin{align}
&\left\{
\begin{aligned}
y&=-\omega^2+O(t+\omega)^{2},\\
z&=-\tfrac 1 3 (1+2\omega^2)\alpha_{2}
	-\tfrac 1 6 (\alpha_{0}-3\alpha_{1}-1)\alpha_{2}(t+\omega)
	+O(t+\omega)^{2},
\end{aligned}\right.\tag{S3-4}\label{sym3-1p}
\\
&\left\{
\begin{aligned} 
y&=-\omega+(1-\alpha_{0})(t+\omega)+O(t+\omega)^2,\\
z&=\tfrac 1 3 (1+2\omega^2)\alpha_{2}+
	\tfrac 1 3(\alpha_{0}-1)\alpha_{2}(t+\omega)+O(t+\omega)^{2},
\end{aligned}\right.\tag{S3-5}\label{sym3-2p}
\\
&\left\{
\begin{aligned}
y&=-\omega+(1+\alpha_{0})(t+\omega)+O(t+\omega)^{2},\\
z&=(t+\omega)^{-1}+\tfrac 1 2 (1+2\omega^2)(\alpha_{0}-\alpha_{1}+1)
	+O(t+\omega),
\end{aligned}\right.\tag{S3-6}\label{sym3-3p}
\end{align}
which correspond to \eqref{sym3-1}, \eqref{sym3-2} and \eqref{sym3-3} 
by $\omega\to\omega^2$.
The actions of $s_{0}$ and $s_{2}$ interchange the symmetric solutions 
as follows: 
\begin{center}
\begin{picture}(190,15)
 \put(0,0){\makebox(10,10){$\scriptstyle s_0$}}
 \put(15,5){\oval(8,8)[l]}
 \put(15,1){\vector(1,0){4}}
 \put(15,9){\vector(1,0){4}}
 \put(20,0){\makebox(30,10){\eqref{sym3-1p}}}
 \put(50,5){\makebox(30,10){$\scriptstyle s_2$}}
 \put(50,5){\vector(1,0){30}}
 \put(80,5){\vector(-1,0){30}}
 \put(80,0){\makebox(30,10){\eqref{sym3-2p}}}
 \put(110,5){\makebox(30,10){$\scriptstyle s_0$}}
 \put(110,5){\vector(1,0){30}}
 \put(140,5){\vector(-1,0){30}}
 \put(140,0){\makebox(30,10){\eqref{sym3-3p}}}
 \put(175,5){\oval(8,8)[r]}
 \put(175,0){\vector(-1,0){4}}
 \put(175,9){\vector(-1,0){4}}
 \put(180,0){\makebox(10,10){$\scriptstyle s_2$}}
\end{picture}
\end{center}
\end{proposition}

\bigskip

In section \ref{MS3}, we will show these symmetric solutions are 
monodromy solvable.  

\subsection{Comparison with classical solutions}
In the case of $\alpha_{0}=0\ (\delta=0)$, the sixth Painlev\'e 
equation \eqref{P6H} admits the Riccati type solution
\begin{align}
\begin{aligned}
y&=t,\\
%=\frac 12 + \left(t-\frac 12\right),\\
t(t-1)\frac{dz}{dt}&=-t(t-1)z^{2}+\Bigl\{1-2t+
\alpha_{3}t+\alpha_{4}(t-1)
\Bigr\}z-(\alpha_{1}+\alpha_{2})\alpha_{2}.
\end{aligned}
\label{Riccati}
\end{align}
If $\alpha_{3}=\alpha_{4}$, \eqref{Riccati} admits a solution
\begin{align*}
y&=1/2 + \left(t-1/2\right),\\
z&=\frac{d}{dt}
\log 
\left[_{2}F_{1}\left(
	\alpha_{2}/2,\,
	(\alpha_{1}+\alpha_{2})/2,\,
	1/2,\,
	4(t-1/2)^2\right) 
\right],
\end{align*}
which is a special case of \eqref{sym2-1}. 
If $\alpha_{1}=\alpha_{3}=\alpha_{4}$, \eqref{Riccati} admits a solution
\begin{align*}
y=&-\omega^2+\left(t+\omega^{2}\right),\\
z=&\frac{d}{dt}
\log
\left[
_{2}F_{1}\left(\frac{1+\alpha_{1}}2,\,
	\frac{1+3\alpha_{1}}6,\,
	\frac 23,\,
	-\left(\frac{-\omega t-1}{t+\omega}\right)^3
\right)
\right]\\
  &\qquad \qquad
  +\frac{t(1-t)(t+\omega^2)+(t^3+3\omega t^{2}+3(1-\omega)t-2)\alpha_{1}}
	{2t(t-1)(t^2-t+1)},
\end{align*}
which is a special case of \eqref{sym3-2}.

\section{Transformation of the Linearization}

Since the Painlev\'e equation is the isomonodromic deformation condition, 
to show the monodromy solvability of a solution of Painlev\'e equation, 
it is sufficient to show the monodromy solvability at a special $t=t_{0}$.
In the following, we will set such special $t=t_{0}$ at the fixed points 
of $\sigma_{1}$ and $\sigma_{2}\circ\sigma_{1}$.

\subsection{Monodromy Solvability of the Symmetric Solution of $\sigma_{1}$}
\label{MS2}
We substitute the symmetric solution \eqref{sym2-1}
into \eqref{OG-1} and take the limit $t\to 1/2$, then
\begin{align}
 \frac{d^2 \psi}{d x^2}+
 \left(
 \frac{1-\alpha_{3}}{x}+\frac{1-\alpha_{3}}{x-1}+\frac{-\alpha_{0}}{x-1/
 2}
\right)\frac{d \psi}{d x}+
\frac{(\alpha_{1}+\alpha_{2})\alpha_{2}}{x(x-1)}
 \psi=0.\label{le2}
\end{align}
This is Heun's equation, whose Riemann scheme is
\[
 P\left\{
\begin{matrix}
0 & 1 & 1/2 & \infty \\
0 & 0 & 0   & \alpha_{2} \\
\alpha_{3} & \alpha_{3} & \alpha_{0}+1 & \alpha_{1}+\alpha_{2}
\end{matrix}; x
\right\}.
\]
Now we set $(2x-1)^2=1-\xi$, then \eqref{le2} becomes
\[
 \frac{d^2 \psi}{d \xi^2}+
 \left(
 \frac{1-\alpha_{3}}{\xi}+\frac{1-\alpha_{0}}{2(\xi-1)}
\right)\frac{d \psi}{d \xi}+
\frac{(\alpha_{1}+\alpha_{2})\alpha_{2}}{4\xi(\xi-1)}
 \psi=0,
\]
which is the hypergeometric equation.

%Then we have the following theorem.
\begin{theorem}
The symmetric solution \eqref{sym2-1} is monodromy 
solvable. Its linearization \eqref{OG-1} at $t=1/2$ is reduced to the 
hyper\-geometric equation: 
\[
 P\left\{
\begin{matrix}
0 & 1 & \infty \\
0 & 0 & \frac{\alpha_{2}}{2}\\
\alpha_{3}&\frac{\alpha_{0}+1}{2}&\frac{\alpha_{1}+\alpha_{2}}{2}
\end{matrix}; \xi
\right\}.
\]
A fundamental solution of \eqref{le2} is given by
\[
 \left({}_{2}F_{1}\left(
	\frac{\alpha_{2}}2, \frac{\alpha_{1}+\alpha_{2}}2, 
	1-\alpha_{3}; 4x(1-x)
\right),\
 (4x(1-x))^{\alpha_{3}}{}_{2}F_{1}\left(
	\frac{\alpha_{2}}2+\alpha_{3},
	\frac{\alpha_{1}+\alpha_{2}}2+\alpha_{3}, 1+\alpha_{3};4x(1-x)
\right)
\right).
\]
\end{theorem}

\bigskip

In the same way, the linearizations of \eqref{sym2-2}, \eqref{sym2-3} 
and \eqref{sym2-4} can be reduced to the hypergeometric equation. 

In section \ref{LM2}, we will calculate the linear monodromy of the 
symmetric solutions of~$\sigma_{1}$ explicitly. 

\subsection{Monodromy Solvability of the Symmetric Solution of 
$\sigma_{2}\circ\sigma_{1}$ }\label{MS3}
We substitute the symmetric solution \eqref{sym3-2} into \eqref{OG-1} 
and take the limit $t\to -\omega^2$, then we have 
\begin{align}
 \frac{d^2 \psi}{d x^2}+
 \left(
	\frac{1-\alpha_{1}}x+\frac{1-\alpha_{1}}{x-1}+
	\frac{-\alpha_{0}}{x+\omega^2}
\right)\frac{d \psi}{d x}+
\left(
\frac{\alpha_{2}(\alpha_{1}+\alpha_{2})}{x(x-1)}-
\frac{2\omega+1}3\frac{\alpha_{0}\alpha_{2}}{x(x-1)(x+\omega^2)}
\right)
 \psi=0.\label{le3}
\end{align}
This is Heun's equation, whose Riemann scheme is
\[
 P\left\{
\begin{matrix}
0 & 1 & -\omega^2 & \infty \\
0 & 0 & 0   & \alpha_{2} \\
\alpha_{1} & \alpha_{1} & 1+\alpha_{0} & \alpha_{1}+\alpha_{2}
\end{matrix}; x
\right\}
=(x+\omega^2)^{-\alpha_{2}}P\left\{
\begin{matrix}
0&1&-\omega^2&\infty\\
0&0&\alpha_{2}&0\\
\alpha_{1}&\alpha_{1}&1+\alpha_{0}+\alpha_{2}&\alpha_{1}
\end{matrix}; x
\right\}.
\]
By the transformation
\[
 \psi=(x+\omega^2)^{-\alpha_{2}}\tilde{\psi},\qquad
 \xi=\omega\frac{x+\omega}{x+\omega^2},\qquad
 \eta=\xi^3,
\]
\eqref{le3} is reduced to the hyper\-geometric equation
\begin{align*}
 \eta(\eta-1)\frac{d^2\tilde{\psi}}{d\eta^2}+
 \left\{\frac 23-\left(1+\frac{1+\alpha_{0}+2\alpha_{2}}3\right)\eta\right\}
	\frac{d\tilde{\psi}}{d\eta}-
 \frac{\alpha_{2}}9\left(
 	1+\alpha_{0}+\alpha_{2}
\right)\tilde{\psi}=0.%\label{le32}
\end{align*}
Then we have the following.

\begin{theorem}
 The symmetric solution \eqref{sym3-2} on $\sigma_{2}\circ\sigma_{1}$ is 
 monodromy solvable. Its linearization 
 \eqref{OG-1} at $t=-\omega^2$ is reduced to the hypergeometric equation  
\[
 P\left\{
\begin{matrix}
0&1&\infty\\
0&0&\alpha_{2}/3\\
1/3&\alpha_{1}&(1+\alpha_{0}+\alpha_{2})/3
\end{matrix};\eta
\right\}.
\]
A fundamental solution is
\[
 \left(
{}_{2}F_{1}\left(\frac{\alpha_{2}}{3},\frac{1}{3}(\alpha_{0}+\alpha_{2}+
1), \frac{2}{3};\eta
	\right),\ \\
\eta^{1/3}{}_{2}F_{1}\left(
\frac{\alpha_{2}+1}{3},\frac 13(\alpha_{0}+\alpha_{2}+2), \frac 43 ; \eta
\right)
 \right).
\]
\end{theorem}

\bigskip

In the same way, \eqref{sym3-1} and \eqref{sym3-3} can be reduced to the 
hypergeometric equation. %Therefore, we can calculate their monodromy. 

In section \ref{LM3}, we will calculate the linear monodromy of the 
symmetric solution of $\sigma_{2}\circ\sigma_{1}$ explicitly.

\section{Linear Monodromy}
%In the following, 
%we take paths $\gamma_{j}$ around singular points $a_{j}$ as a standard 
%way in the figure.
%\begin{center}
%\includegraphics{path.eps}
%\end{center}
%We represent a choice of paths corresponding to the monodromy matrices 
%$M_{j}$ around the singular points $a_{j}$ as
%\[
% \left[\begin{matrix}
%a_{1} & a_{2} & \cdots & a_{k}\\
%M_{1} & M_{2} & \cdots & M_{k}
%\end{matrix}
%\right],
%\]
%where the fundamental solution $\Psi=(\psi_{1},\psi_{2})$ is 
%analytically continued along  
%$\gamma_{j}$ as 
%\[
% \Psi(\gamma_{j}(x))=\Psi(x)M_{j}.
%\]
%The monodromy matrices satisfy the following relation
%\[
% M_{k}\cdots M_{2}M_{1}=1.
%\]

\subsection{Linear Monodromy of the Symmetric Solution of $\sigma_{1}$}
\label{LM2}

By $(2x-1)^2=1-\xi$, 
$\xi$-surface $\mathbb{P}_{1}\setminus\{0,1,\infty\}$ is double covered 
by the $x$-surface 
$\mathbb{P}_{1}\setminus\{0,1/2,1,\infty\}$.
Let us make a cut from $\xi=1$ to $\infty$ ($x=1/2$ to $\infty$), and 
suppose that $x=1$ is on the second $\xi$-plane, and the base point and 
$x=0$ are on the first $\xi$-plane. 

Let us denote the standard paths around 
$x=0,1/2,1,\infty$ by
$\gamma_{0},\gamma_{1/2},\gamma_{1},\gamma_{\infty}$,
and denote the standard paths around $\xi=0,1,\infty$ by
$L_{0},L_{1},L_{\infty}$.

\begin{center}
\includegraphics[width=12cm]{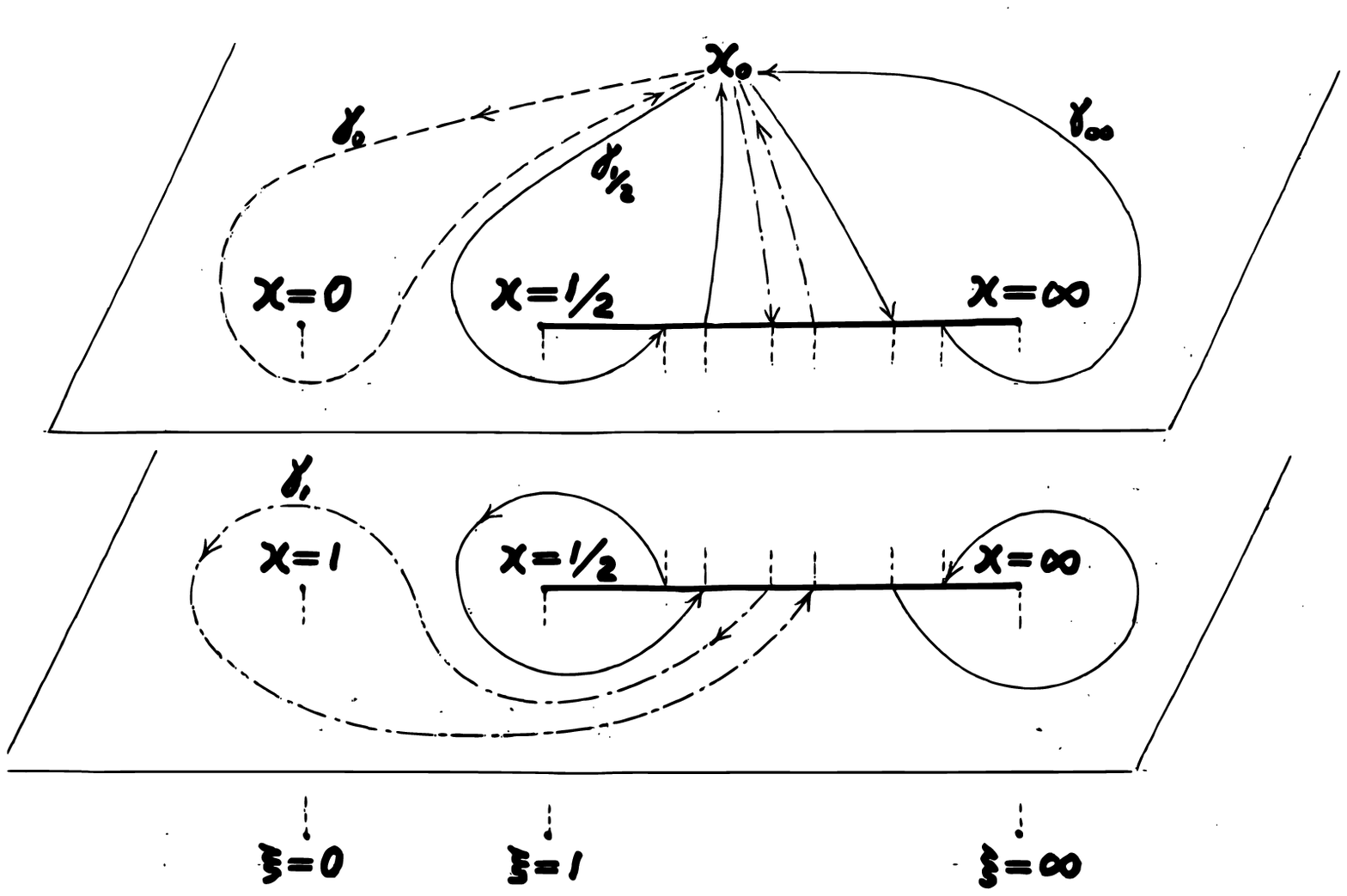}\\
\includegraphics[width=12cm]{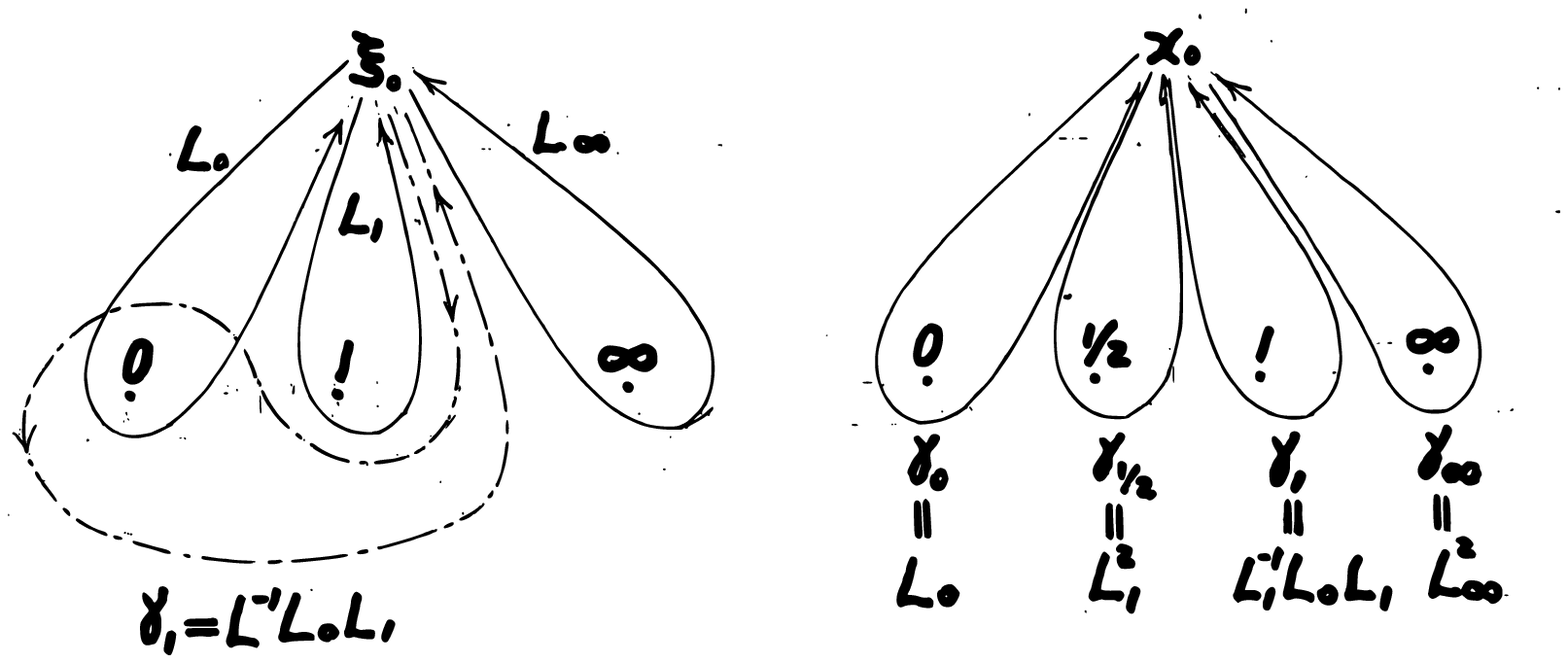}
\end{center}

Then we have 
$$
\gamma_{0}=L_{0}, \qquad
\gamma_{1/2}=L_{1}^2, \qquad
\gamma_{1}=L_{1}^{-1}L_{0}L_{1},  \qquad
\gamma_{\infty}=L_{\infty}^2,
$$
which satisfy
\[
 \gamma_{0}\gamma_{1/2}\gamma_{1}\gamma_{\infty}=
L_{0}(L_{1}^2)(L_{1}^{-1}L_{0}L_{1})L_{\infty}^2=1. 
\]
%Therefore,
%for a fundamental solution 
%$\Psi=(\psi_{1},\psi_{2})$ of \eqref{OG-1},
%their monodromy is represented by
%$$
%M_{0}=G_{0},\qquad
%M_{t}=G_{1}^2, \qquad
%M_{1}=G_{1}G_{0}G_{1}^{-1}, \qquad
%M_{\infty}=G_{\infty}^2,
%$$
%which satisfies
%\[
% M_{\infty}M_{1}M_{t}M_{0}=
%G_{\infty}^{2}(G_{1}G_{0}G_{1}^{-1})G_{1}^2G_{0}=1.
%\]

Then we obtain monodromy matrix $M_{j}$ along path $\gamma_{j}$ 
($j=0, 1/2, 1, \infty$). 
\begin{theorem}\label{mon2-thm}
The linear monodromy of the symmetric solution \eqref{sym2} of $\sigma_{1}$
is represented as
\begin{gather*}
 M_{0}=\Gamma_{{0\infty}}\Lambda_{0}\Gamma_{0\infty}^{-1},\qquad
 M_{1/2}=\Gamma_{{\frac 12\infty}}
	\Lambda_{\frac 12}^2\Gamma_{\frac 12\infty}^{-1},\\
 M_{1}=(\Gamma_{{\frac 12\infty}}\Lambda_{\frac 12}\Gamma_{\frac 12\infty}^{-1})
	(\Gamma_{{0\infty}}\Lambda_{1}\Gamma_{0\infty}^{-1})
	(\Gamma_{{\frac 12\infty}}
\Lambda_{\frac 12}^{-1}\Gamma_{\frac 12\infty}^{-1}),\qquad
 M_{\infty}=(e^{2\pi i T_{\infty}})^2,
\end{gather*}
where
\begin{gather*}
\Lambda_{0}=\Lambda_{1}=\begin{pmatrix}
1&0\\0&e^{2\pi i \alpha_{3}}
\end{pmatrix},\qquad
\Lambda_{\frac 12}=\begin{pmatrix}
1&0\\
0&-e^{\pi i \alpha_{0}}
\end{pmatrix},\qquad
e^{2\pi i T_{\infty}}=\begin{pmatrix}
e^{\pi i\alpha_{2}}&0\\0&e^{\pi i (\alpha_{1}+\alpha_{2})}
\end{pmatrix},\\
\Gamma_{0\infty}=\Gamma_{1\infty}=\begin{pmatrix}
e^{-\alpha_{2}\pi i/2}
\frac{\Gamma(1-\alpha_{3})\Gamma(\alpha_{1}/2)}
{\Gamma((\alpha_{1}+\alpha_{2})/2)\Gamma(1-\alpha_3-\alpha_{2}/2)}
&	%%%%%%%%%%%%%%%%%%%%%%%%%
e^{-(\alpha_{3}+\alpha_{2}/2)\pi i}\frac
{\Gamma(1+\alpha_{3})\Gamma(\alpha_{1}/2)}
{\Gamma(1-\alpha_{2}/2)\Gamma((1-\alpha_{0}-\alpha_{2})/2)}
\\[+5pt]	%%%%%%%%%%%%%%%%%%%%%%%%
e^{-(\alpha_{1}+\alpha_{2})\pi i/2}\frac
{\Gamma(1-\alpha_{3})\Gamma(-\alpha_{1}/2)}
{\Gamma(\alpha_{2}/2)\Gamma((1+\alpha_{0}+\alpha_{2})/2)}
& %%%%%%%%%%%%%%%%%%%%%%%%%%%
e^{(\alpha_{0}+\alpha_{2}-1)\pi i/2}\frac
{\Gamma(1+\alpha_{3})\Gamma(-\alpha_{1}/2)}
{\Gamma((1-\alpha_{1}-\alpha_{2})/2)\Gamma(\alpha_{3}+\alpha_{2}/2)}
\end{pmatrix},
\\ %%%%%%%%%%%%%%%%%%%%%%%%%%%%%
\Gamma_{\frac 12 \infty}=\begin{pmatrix}
\frac{\Gamma((3+\alpha_{0})/2)\Gamma(-\alpha_{1}/2)}
{\Gamma((1+\alpha_{0}+\alpha_{2})/2)\Gamma((2-\alpha_{1}-\alpha_{2})/2)}
& %%%%%%%%%%%%%%%%%%%%%%%%%%%%%%
-\frac{\Gamma((3+\alpha_{0})/2)\Gamma(\alpha_{1}/2)}
{\Gamma(1-\alpha_{3}-\alpha_{2}/2)\Gamma(1-\alpha_{2}/2)}
\\[+5pt] %%%%%%%%%%%%%%%%%%%%%%%%%
-e^{\pi i ((1+\alpha_{0})/2)}\frac
{\Gamma((1-\alpha_{0})/2)\Gamma(\alpha_{1}/2)}
{\Gamma(\alpha_{2}/2)\Gamma(\alpha_{3}+\alpha_{2}/2)}
& %%%%%%%%%%%%%
e^{\pi i ((1+\alpha_{0})/2)}\frac
{\Gamma((1-\alpha_{0})/2)\Gamma(\alpha_{1}/2)}
{\Gamma((\alpha_{1}+\alpha_{2})/2)\Gamma((1-\alpha_{0}-\alpha_{2})/2)}
\end{pmatrix}.
\end{gather*}
\end{theorem}

\subsection{Linear Monodromy of the Symmetric Solution of 
$\sigma_{2}\circ\sigma_{1}$} \label{LM3}

Let us denote the standard paths around $\xi=\omega, \omega^2, 1, 
\infty$ by $\gamma_{\omega}$, $\gamma_{\omega^2}$, $\gamma_{1}$, 
$\gamma_{\infty}$, and the standard paths around $\eta=0, 1, \infty$ by 
$L_{0}$, $L_{1}$, $L_{\infty}$.% We set on $x$-space
%\[
% \begin{bmatrix}
%0 & t & 1 & \infty\\
%M_{0}&M_{t}&M_{1}&M_{\infty}
%\end{bmatrix},
%\]
%on $\xi$-space
%\[
% \begin{bmatrix}
% \omega & \omega^2 & 1 & \infty\\
% \tilde{M}_{\omega} & \tilde{M}_{\omega^2} & \tilde{M}_{1} & \tilde{M}_{\infty}
%\end{bmatrix},
%\]
%and on $\eta$-space
%\[
% \begin{bmatrix}
%0&1&\infty\\
%G_{0}&G_{1}&G_{\infty}
%\end{bmatrix}.
%\]
%Then 
%\[
% M_{0}=\tilde{M}_{1}
%\]
%%\begin{center}
%%\includegraphics{xi30.eps} \hspace{10mm}
%%\includegraphics{eta30.eps}
%%\end{center}

The $\xi$-space and the $\eta$-space is connected by $\xi^3=\eta$. 

\begin{center}
\includegraphics{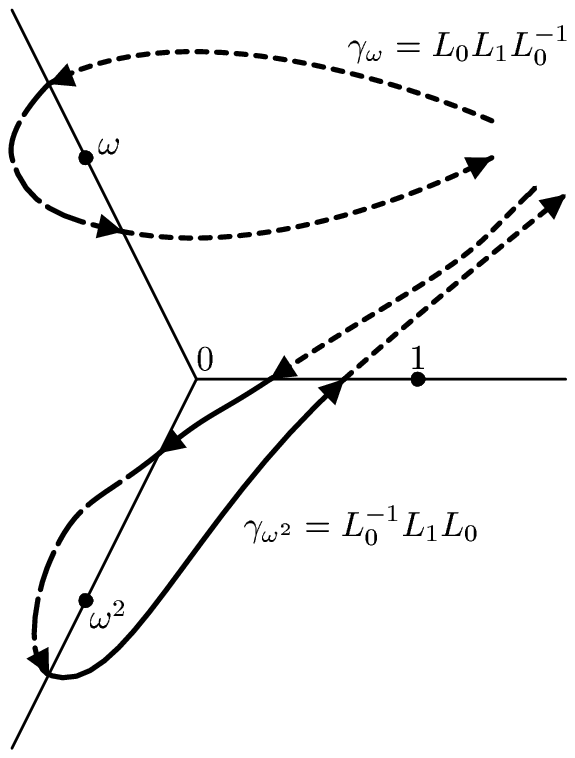}\hspace{10mm}
\includegraphics{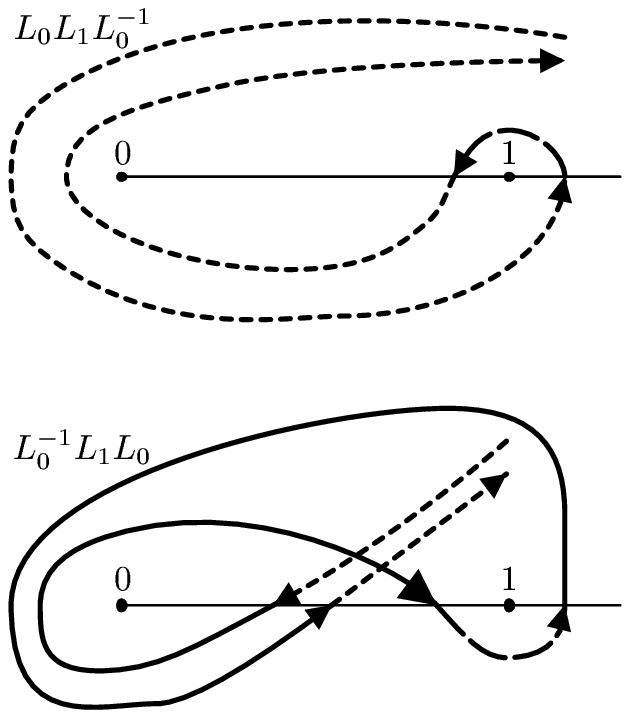}
\end{center}

\begin{center}
\includegraphics{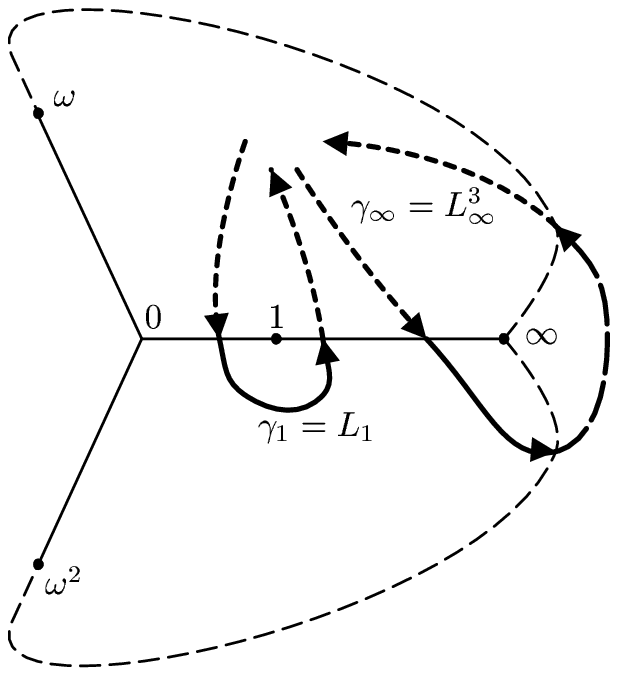}\hspace{10mm}
\includegraphics{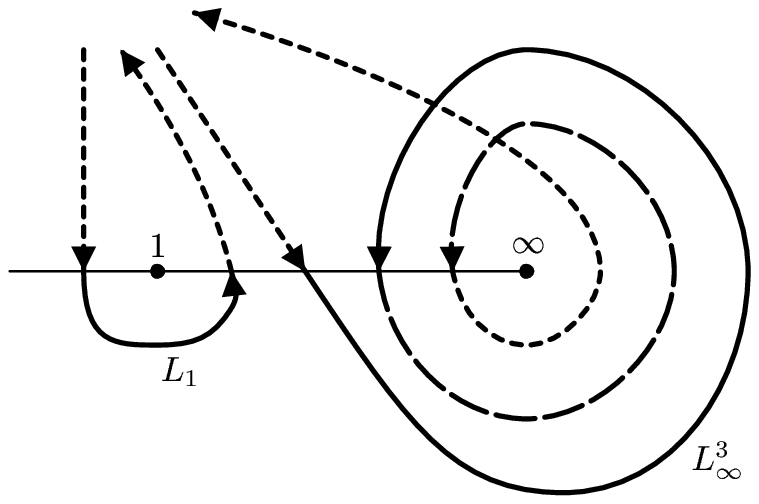}
\end{center}

Then we have 
\[
 \gamma_{\omega}=L_{0}L_{1}L_{0}^{-1},\qquad
 \gamma_{\omega^{2}}=L_{0}^{-1}L_{1}L_{0},\qquad
 \gamma_{1}=L_{1},\qquad
 \gamma_{\infty}=L_{\infty}^3,
\]
where 
\[
 \gamma_{\omega}\gamma_{\omega^2}\gamma_{1}\gamma_{\infty}=
 (L_{0}L_{1}L_{0}^{-1})(L_{0}^{-1}L_{1}L_{0})L_{1}(L_{\infty}^3)
 =L_{0}L_{1}L_{0}^{-3}L_{\infty}.
\]
Since the monodromy along the path $L_{0}^{3}$ is $1$, the monodromy 
along the path 
$\gamma_{\omega}\gamma_{\omega^2}\gamma_{1}\gamma_{\infty}$ is $1$.

Then we obtain monodromy matrix $M_{j}$ along path $\gamma_{j}$ 
($j=\omega, \omega^2, 1, \infty$). 

\begin{theorem} \label{mon3-thm}
The linear monodromy of the symmetric solution \eqref{sym3-2} of $\sigma_{2}\circ\sigma_{1}$
is represented as follows:
\begin{gather*}
 M_{\omega^2}=e^{-2\pi i \alpha_{2}}\left(
	e^{2\pi i T_{\infty}}
\right)^{3},\qquad
 M_{0}=(\Gamma_{1\infty}\Lambda_{1}\Gamma_{1\infty}^{-1}),
\\
M_{1}=(\Gamma_{0\infty}\Lambda_{0}\Gamma_{0\infty}^{-1})
(\Gamma_{1\infty}\Lambda_{1}\Gamma_{1\infty}^{-1})
(\Gamma_{0\infty}\Lambda_{0}^{-1}\Gamma_{0\infty}^{-1}),\\
M_{\infty}=e^{2\pi i\alpha_{2}}
(\Gamma_{0\infty}\Lambda_{0}^{-1}\Gamma_{0\infty}^{-1})
(\Gamma_{1\infty}\Lambda_{1}\Gamma_{1\infty}^{-1})
(\Gamma_{0\infty}\Lambda_{0}\Gamma_{0\infty}^{-1}),
\end{gather*}
where
\begin{gather*}
\Lambda_{0}=\begin{pmatrix}1&0\\0&e^{2\pi i/3}\end{pmatrix},\qquad
\Lambda_{1}=\begin{pmatrix}1&0\\0&e^{2\pi i
	(1-\alpha_{0}-2\alpha_{2})/3}\end{pmatrix},\qquad
e^{2\pi i T_{\infty}}=\begin{pmatrix}
e^{2\pi i\alpha_{2}/3}&0\\0&e^{2\pi i(1+\alpha_{0}+\alpha_{2})/3}
\end{pmatrix},\\ %%%%%%%%%%%%%%%%%%%%%%
\Gamma_{0\infty}=\begin{pmatrix}
e^{\alpha_{2}\pi i/3}\frac{\Gamma(2/3)\Gamma((1+\alpha_{0})/3)}
{\Gamma((1+\alpha_{0}+\alpha_{2})/3)\Gamma((2-\alpha_{2})/3)}
& %%%%%%%%%%%%%%%%%%%%%%%%%%%%%%%%%%%%%%%%%%%%%%%%%%%%
e^{(1+\alpha_{2})\pi i/3}\frac
{\Gamma(4/3)\Gamma((1+\alpha_{0})/3)}
{\Gamma(1-\alpha_{2}/3)\Gamma((2+\alpha_{0}+\alpha_{2})/3)}
\\[+5pt] %%%%%%%%%%%%%%%%%%%%%%%%%%%%%%%%%%%%%%%%%%%%%%%%%%%
e^{(1+\alpha_{0}+\alpha_{2})\pi i/3}\frac
{\Gamma(2/3)\Gamma((1+\alpha_{0})/(-3))}
{\Gamma(\alpha_{2}/3)\Gamma((1-\alpha_{0}-\alpha_{2})/3)}
& %%%%%%%%%%%%%%%%%%%%%%%%%%%%%%%%%%%%%%%%%%%%%%%%%%%%
e^{(2+\alpha_{0}+\alpha_{2})\pi i/3}\frac
{\Gamma(4/3)\Gamma((1+\alpha_{0})/(-3))}
{\Gamma((2-\alpha_{0}-\alpha_{2})/3)\Gamma((1+\alpha_{2})/3)}
\end{pmatrix},\\ %%%%%%%%%%%%%%%%%%%%%%%%%%%%%%%%%%%%%%%
\Gamma_{1\infty}=\begin{pmatrix}
\frac{\Gamma((2+\alpha_{0}+2\alpha_{2})/3)\Gamma((1+\alpha_{0})/3)}
{\Gamma((1+\alpha_{0}+\alpha_{2})/3)\Gamma((2+\alpha_{0}+\alpha_{2})/3)}
& %%%%%%%%%%%%%%%%%%%%%%%%%%%%%%%%%%%%%%%%
e^{(-1+\alpha_{0}+2\alpha_{2})\pi i/3}
\frac{\Gamma((4-\alpha_{0}-2\alpha_{2})/3)\Gamma((1+\alpha_{0})/3)}
{\Gamma((2-\alpha_{2})/3)\Gamma(1-\alpha_{2}/3)}
\\[+5pt] %%%%%%%%%%%%%%%%%%%%%%%%%%%%%%%%%%%%%%%%%%%%%%%%
\frac{\Gamma((2+\alpha_{0}+2\alpha_{2})/3)\Gamma((1+\alpha_{0})/(-3))}
{\Gamma(\alpha_{2}/3)\Gamma((1+\alpha_{2})/3)}
& %%%%%%%%%%%%%%%%%%%%%%%%%%%%%%%%%%
e^{(-1+\alpha_{0}+2\alpha_{2})\pi i/3}\frac
{\Gamma((4-\alpha_{0}-2\alpha_{2})/3)\Gamma((1+\alpha_{0})/(-3))}
{\Gamma((1-\alpha_{0}-\alpha_{2})/3)\Gamma((2-\alpha_{0}-\alpha_{2})/3)}
\end{pmatrix}.
\end{gather*}
\end{theorem}

\section{Characterizations of the monodromy}\label{char}

In this section, we will characterize the linear monodromy of our 
symmetric solutions on Fricke's cubic surface of monodromy. 

\bigskip

To normalize the monodromy matrices in $SL(2,\mathbb{C})$, 
we replace $\psi$ of the 
linearization \eqref{OG-1} as 
\[
 \psi \to x^{\alpha_{4}/2}(x-1)^{\alpha_{3}/2}(x-t)^{(\alpha_{0}-1)/2}\psi,
\]
then the Riemann scheme of \eqref{OG-1} becomes
\[
P\left\{x; 
 \begin{matrix}
 0	&1	&t	&y	&\infty\\
-\frac{\alpha_{4}}2&-\frac{\alpha_{3}}2&-\frac{\alpha_{0}-1}2&0&
		-\frac{\alpha_{1}}2\\ 
\frac{\alpha_{4}}2&\frac{\alpha_{3}}2&\frac{\alpha_{0}-1}2+1&2&
		\frac{\alpha_{1}}2
\end{matrix}\right\},
\]
and then the monodromy matrices 
$\left\{M_{0},M_{t},M_{1},M_{\infty}\right\}$
become elements of 
$SL(2,\mathbb{C})$. Hereafter we use this linearization.

\subsection{Linear monodromy of the Symmetric Solutions of $\sigma_{1}$}
\label{char2}

From the theorem \ref{mon2-thm}, we can calculate 
$\{p_{0}, p_{1}, p_{t}, p_{\infty}; p_{01}, p_{1t}, p_{t0}\}$ for the 
symmetric solution \eqref{sym2-1}. We have
\[
 p_{0}=p_{1}=2\cos\pi\alpha_{3}, \qquad
 p_{t}=2\cos\pi(\alpha_{0}+1), \qquad
 p_{\infty}=2\cos\pi\alpha_{1},
\]
 and
\begin{align}
&p_{1t}=p_{t0}=2\left(
\cos\pi(\alpha_{3}+1)+\cos\frac{\pi}2(\alpha_{0}-\alpha_{1}-1)
	+\cos\frac{\pi}2(\alpha_{0}+\alpha_{1}-1)
\right),\label{X21}\\
&p_{01}=-2\left(1+\cos\pi(\alpha_{0}+1)+\cos\pi\alpha_{1}+
4\cos\frac\pi 2\alpha_{1}\cos\pi\alpha_{3}\cos\frac\pi 2(\alpha_{0}+1)\right).
\end{align}
In the same way, for the solution \eqref{sym2-2}, we have
\begin{align}
& p_{0}=p_{1}=2\cos\pi\alpha_{3}, \qquad
 p_{t}=2\cos\pi(\alpha_{0}+1), \qquad
 p_{\infty}=2\cos\pi\alpha_{1},\notag \\
&p_{1t}=p_{t0}=2\left(
\cos\pi(\alpha_{3}+1)+\cos\frac{\pi}2(\alpha_{0}-\alpha_{1}+1)
	+\cos\frac{\pi}2(\alpha_{0}+\alpha_{1}+1)
\right),\label{X22}\\
&p_{01}=-2\left(1+\cos\pi(\alpha_{0}-1)+\cos\pi\alpha_{1}+
4\cos\frac\pi 2\alpha_{1}\cos\pi\alpha_{3}\cos\frac\pi 2(\alpha_{0}-1)\right),
\end{align}
for \eqref{sym2-3}, we have
\begin{align}
& p_{0}=p_{1}=2\cos\pi\alpha_{3}, \qquad
 p_{t}=2\cos\pi(\alpha_{0}+1), \qquad
 p_{\infty}=2\cos\pi\alpha_{1},\notag \\
&p_{1t}=p_{t0}=2\left(
\cos\pi\alpha_{3}+\cos\frac{\pi}2(\alpha_{0}-\alpha_{1}-1)
	+\cos\frac{\pi}2(\alpha_{0}+\alpha_{1}+1)
\right),\label{X23}\\
&p_{01}=-2\left(1+\cos\pi\alpha_{0}+\cos\pi(\alpha_{1}+1)+
4\cos\frac\pi 2(\alpha_{1}+1)\cos\pi\alpha_{3}\cos\frac\pi 2\alpha_{0}\right),
\end{align}
and, for \eqref{sym2-4}, we have
\begin{align}
& p_{0}=p_{1}=2\cos\pi\alpha_{3}, \qquad
 p_{t}=2\cos\pi(\alpha_{0}+1), \qquad
 p_{\infty}=2\cos\pi\alpha_{1},\notag \\
&p_{1t}=p_{t0}=2\left(
\cos\pi\alpha_{3}+\cos\frac{\pi}2(\alpha_{0}-\alpha_{1}+1)
	+\cos\frac{\pi}2(\alpha_{0}+\alpha_{1}-1)
\right),\label{X24}\\
&p_{01}=-2\left(1+\cos\pi\alpha_{0}+\cos\pi(\alpha_{1}-1)+
4\cos\frac\pi 2(\alpha_{1}-1)\cos\pi\alpha_{3}\cos\frac\pi 2\alpha_{0}\right).
\end{align}
Certainly, they satisfies \eqref{eq-mon}. 

\bigskip

Conversely, if we set 
\begin{gather*}
X= p_{1t}=p_{t0}, \qquad \qquad Y= p_{01},\\
 A=p_{0}=p_{1}=2\cos\pi\alpha_{3}, \qquad
 B=p_{t}=2\cos\pi(\alpha_{0}+1), \qquad
 C=p_{\infty}=2\cos\pi\alpha_{1},
\end{gather*}
then \eqref{eq-mon} becomes
\[
 X^2Y+Y^2+2X^2-(A^2+BC)Y-2(AB+AC)X+2A^2+B^2+C^2+A^2BC-4=0.
\]
As an equation of $Y$, the discriminant of this equation can be 
factorized as follows 
\[
 (X-X_{1})(X-X_{2})(X-X_{3})(X-X_{4}),
\]
where $X_{1}$, $X_{2}$, $X_{3}$ and $X_{4}$ are the right hand sides of 
\eqref{X21}, \eqref{X22}, \eqref{X23} and \eqref{X24}. 

Now, we have the following.
\begin{theorem}
Under the constraints $p_{0}=p_{1}$ and $p_{0t}=p_{t1}$, the relation of 
monodromy \eqref{eq-mon} admits a double root as an equation of $p_{01}$, 
if and only if 
$\left\{M_{0}, M_{t}, M_{1}, M_{\infty}\right\}$ is the linear monodromy 
of a symmetric solution of $\sigma_{1}$.
\end{theorem}

\subsection{Linear monodromy of the Symmetric Solutions of $\sigma_{2}\circ\sigma_{1}$}
\label{char3}

From the theorem \ref{mon3-thm}, we can calculate 
$\{p_{0}, p_{1}, p_{t}, p_{\infty}; p_{01}, p_{1t}, p_{t0}\}$ for the 
symmetric solution \eqref{sym3-2}. We have
\begin{align}
& p_{0}=p_{1}=p_{\infty}=2\cos(\pi\alpha_{1}),\qquad 
 p_{t}=2\cos\pi(\alpha_{0}-1),\\
&p_{0t}=p_{1t}=p_{t\infty}=-1-2\cos\left(\tfrac{2\pi}3\alpha_{0}\right)
	+4\cos\left(\tfrac\pi 3 \alpha_{0}\right)\cos\pi\alpha_{1},
\label{X31}
\end{align}
In the same way, for the solution \eqref{sym3-1}, we have
\begin{align}
& p_{0}=p_{1}=p_{\infty}=2\cos(\pi\alpha_{1}),\qquad 
 p_{t}=2\cos\pi(\alpha_{0}-1),\\
&p_{0t}=p_{1t}=p_{t\infty}=-1-2\cos\left(\tfrac{2\pi}3(\alpha_{0}-2)\right)
	+4\cos\left(\tfrac\pi 3 (\alpha_{0}-2)\right)\cos\pi\alpha_{1},
\label{X32}
\end{align}
and, for \eqref{sym3-3}, we have
\begin{align}
& p_{0}=p_{1}=p_{\infty}=2\cos(\pi\alpha_{1}),\qquad 
 p_{t}=2\cos\pi(\alpha_{0}-1),\\
&p_{0t}=p_{1t}=p_{t\infty}=-1-2\cos\left(\tfrac{2\pi}3(\alpha_{0}+2)\right)
	+4\cos\left(\tfrac\pi 3 (\alpha_{0}+2)\right)\cos\pi\alpha_{1}.
\label{X33}
\end{align}

Now, 
$\{p_{0}, p_{1}, p_{t}, p_{\infty}; p_{\infty 1}, p_{t1}, p_{t \infty}\}$
also parameterize the monodromy, and they satisfy the following relation
\begin{multline}
p_{\infty 1}p_{1t}p_{t\infty}+p_{\infty 1}^2+p_{1 t}^2+p_{t\infty}^2
-(p_{\infty}p_{1}+p_{t}p_{0})
p_{\infty 
1}-(p_{1}p_{t}+p_{\infty}p_{0})p_{1t}-(p_{t}p_{\infty}+p_{1}p_{0})p_
{t\infty} 
\\
+p_{\infty}^2+p_{1}^2+p_{t}^2+p_{0}^{2}
+p_{0}p_{t}p_{1}p_{t}p_{\infty}-4=0.
\label{eq-mon2}
\end{multline}
Since $M_{0}$, $M_{t}$, $M_{1}$, $M_{\infty}$ are elements of 
$SL(2,\mathbb{C})$, we have 
\[
 p_{0t}=\mathrm{tr}M_{0}M_{t}=\mathrm{tr}\left(M_{1}M_{\infty}\right)^{-1}=
\mathrm{tr}M_{\infty}M_{1}=p_{\infty 1}.
\]
And thus, 
$\{p_{0}, p_{1}, p_{t}, p_{\infty}; p_{t 0}, p_{t 1}, p_{t \infty}\}$
also parameterize the monodromy. 

\bigskip

Conversely, if we set
\[
 X=p_{t0}=p_{t1}=p_{t\infty}, \qquad 
 A=p_{0}=p_{1}=p_{\infty}=2\cos(\pi\alpha_{1}), \qquad 
 B= p_{t}=2\cos\pi(\alpha_{0}-1),
\]
then \eqref{eq-mon2} becomes
\[
 X^3+3X^2-3(A^2+AB)X+3A^2+B^2+A^3B-4=0.
\]
We can factorize the left hand side of this equation as follows
\[
 (X-X_{1})(X-X_{2})(X-X_{3}),
\]
where $X_{1}$, $X_{2}$, $X_{3}$ are the right hand sides of 
\eqref{X31}, \eqref{X32}, \eqref{X33}.

Now we have the following theorem. 
\begin{theorem}
Under the constraints $p_{0}=p_{1}=p_{\infty}$ and 
$p_{t0}=p_{t1}=p_{t\infty}$, the relation \eqref{eq-mon2} becomes a 
third order equation of $p_{t0}$, whose solutions correspond to 
the linear monodromy of the symmetric solutions of 
$\sigma_{2}\circ\sigma_{1}$. 
\end{theorem}

\end{document}